\documentclass[11pt]{article} 
\usepackage{latexsym} 
\usepackage{amsmath}
\usepackage{amssymb} 
\usepackage{amsfonts}
\usepackage{amscd} 
\usepackage{graphicx}
\usepackage{layout}
\begin{document} 
\newtheorem{Th}{Theorem}[section]
\newtheorem{Cor}{Corollary}[section]
\newtheorem{Prop}{Proposition}[section]
\newtheorem{Lem}{Lemma}[section]
\newtheorem{Def}{Definition}[section]
\newtheorem{Rem}{Remark}[section]
\newtheorem{Ex}{Example}[section]
\newtheorem{stw}{Proposition}[section]


\newcommand{\bet}{\begin{Th}}
\newcommand{\ent}{\stepcounter{Cor}
   \stepcounter{Prop}\stepcounter{Lem}\stepcounter{Def}
   \stepcounter{Rem}\stepcounter{Ex}\end{Th}}


\newcommand{\bec}{\begin{Cor}}
\newcommand{\enc}{\stepcounter{Th}
   \stepcounter{Prop}\stepcounter{Lem}\stepcounter{Def}
   \stepcounter{Rem}\stepcounter{Ex}\end{Cor}}
\newcommand{\bep}{\begin{Prop}}
\newcommand{\enp}{\stepcounter{Th}
   \stepcounter{Cor}\stepcounter{Lem}\stepcounter{Def}
   \stepcounter{Rem}\stepcounter{Ex}\end{Prop}}
\newcommand{\bel}{\begin{Lem}}
\newcommand{\enl}{\stepcounter{Th}
   \stepcounter{Cor}\stepcounter{Prop}\stepcounter{Def}
   \stepcounter{Rem}\stepcounter{Ex}\end{Lem}}
\newcommand{\bef}{\begin{Def}}
\newcommand{\enf}{\stepcounter{Th}
   \stepcounter{Cor}\stepcounter{Prop}\stepcounter{Lem}
   \stepcounter{Rem}\stepcounter{Ex}\end{Def}}
\newcommand{\ber}{\begin{Rem}}
\newcommand{\enr}{
   \stepcounter{Th}\stepcounter{Cor}\stepcounter{Prop}
   \stepcounter{Lem}\stepcounter{Def}\stepcounter{Ex}\end{Rem}}
\newcommand{\bee}{\begin{Ex}}
\newcommand{\ene}{
   \stepcounter{Th}\stepcounter{Cor}\stepcounter{Prop}
   \stepcounter{Lem}\stepcounter{Def}\stepcounter{Rem}\end{Ex}}
\newcommand{\Proof}{\noindent{\it Proof\,}:\ }
\newcommand{\beP}{\Proof}
\newcommand{\enP}{\hfill $\Box$ \par\vspace{5truemm}}

\newcommand{\EE}{\mathbf{E}}
\newcommand{\QQ}{\mathbf{Q}}
\newcommand{\R}{\mathbf{R}}
\newcommand{\C}{\mathbf{C}}
\newcommand{\ZZ}{\mathbf{Z}}
\newcommand{\KK}{\mathbf{K}}
\newcommand{\NN}{\mathbf{N}}
\newcommand{\PP}{\mathbf{P}}
\newcommand{\HH}{\mathbf{H}}
\newcommand{\uuu}{\boldsymbol{u}}
\newcommand{\xxx}{\boldsymbol{x}}
\newcommand{\aaa}{\boldsymbol{a}}
\newcommand{\bbb}{\boldsymbol{b}}
\newcommand{\AAA}{\mathbf{A}}
\newcommand{\BBB}{\mathbf{B}}
\newcommand{\ccc}{\boldsymbol{c}}
\newcommand{\iii}{\boldsymbol{i}}
\newcommand{\jjj}{\boldsymbol{j}}
\newcommand{\kkk}{\boldsymbol{k}}
\newcommand{\rrr}{\boldsymbol{r}}
\newcommand{\FFF}{\boldsymbol{F}}
\newcommand{\yyy}{\boldsymbol{y}}
\newcommand{\ppp}{\boldsymbol{p}}
\newcommand{\qqq}{\boldsymbol{q}}
\newcommand{\nnn}{\boldsymbol{n}}
\newcommand{\vvv}{\boldsymbol{v}}
\newcommand{\eee}{\boldsymbol{e}}
\newcommand{\fff}{\boldsymbol{f}}
\newcommand{\www}{\boldsymbol{w}}
\newcommand{\0}{\boldsymbol{0}}
\newcommand{\lon}{\longrightarrow}
\newcommand{\ga}{\gamma}
\newcommand{\pa}{\partial}
\newcommand{\QED}{\hfill $\Box$}
\newcommand{\id}{{\mbox {\rm id}}}
\newcommand{\Ker}{{\mbox {\rm Ker}}}
\newcommand{\grad}{{\mbox {\rm grad}}}
\newcommand{\ind}{{\mbox {\rm ind}}}
\newcommand{\rot}{{\mbox {\rm rot}}}
\newcommand{\diver}{{\mbox {\rm div}}}
\newcommand{\Gr}{{\mbox {\rm Gr}}}
\newcommand{\LG}{{\mbox {\rm LG}}}
\newcommand{\Diff}{{\mbox {\rm Diff}}}
\newcommand{\Symp}{{\mbox {\rm Symp}}}
\newcommand{\Ct}{{\mbox {\rm Ct}}}
\newcommand{\Uns}{{\mbox {\rm Uns}}}
\newcommand{\rank}{{\mbox {\rm rank}}}
\newcommand{\sign}{{\mbox {\rm sign}}}
\newcommand{\Spin}{{\mbox {\rm Spin}}}
\newcommand{\Sp}{{\mbox {\rm sp}}}
\newcommand{\Int}{{\mbox {\rm Int}}}
\newcommand{\Hom}{{\mbox {\rm Hom}}}
\newcommand{\codim}{{\mbox {\rm codim}}}
\newcommand{\ord}{{\mbox {\rm ord}}}
\newcommand{\Iso}{{\mbox {\rm Iso}}}
\newcommand{\corank}{{\mbox {\rm corank}}}
\def\mod{{\mbox {\rm mod}}}
\newcommand{\pt}{{\mbox {\rm pt}}}
\newcommand{\qed}{\hfill $\Box$ \par}
\newcommand{\spe}{\vspace{0.4truecm}}
\newcommand{\ad}{{\mbox{\rm ad}}}
%
\newenvironment{FRAME}{\begin{trivlist}\item[]
	\hrule
	\hbox to \linewidth\bgroup
		\advance\linewidth by -10pt
		\hsize=\linewidth
		\vrule\hfill
		\vbox\bgroup
			\vskip5pt
			\def\thempfootnote{\arabic{mpfootnote}}
			\begin{minipage}{\linewidth}}{%
			\end{minipage}\vskip5pt
		\egroup\hfill\vrule
	\egroup\hrule
	\end{trivlist}}

\title{
Duality on geodesics of Cartan distributions 
\\
and sub-Riemannian pseudo-product structures
} 

\author{G. Ishikawa\thanks{This work was supported by KAKENHI No.22340030 and No.23654058.}, Y. Kitagawa and W. Yukuno}


\date{ }

\maketitle

\begin{abstract} 
Given a five dimensional space endowed with a Cartan distribution, the abnormal geodesics 
form another five dimensional space with a cone structure. Then it is shown in \cite{IKY}, 
that, if the cone structure is regarded as a control system, then, 
the space of abnormal geodesics of the cone structure is naturally identified with the original space. 
In this paper, we provide an exposition on the duality by abnormal geodesics in a wider framework, 
namely, in terms of quotients of control systems and sub-Riemannian pseudo-product structures. 
Also we consider the controllability of cone structures and describe the constrained Hamiltonian equations 
on normal and abnormal geodesics. 
\end{abstract}

\section{Introduction.}

A distribution on a five dimensional manifold is called a {\it Cartan distribution} if its growth is $(2, 3, 5)$ (see \S \ref{Cartan distributions}). 
Cartan distributions were studied by Cartan in the famous paper \cite{Cartan}. 
Cartan distributions arise in many problems and they are studied, after E. Cartan, 
in detail from various viewpoints by many mathematicians (see for instance, 
Yamaguchi\cite{Yamaguchi}, Bryant and Hsu \cite{BH}, Zhitomirskii\cite{Zhitomirskii}, 
Nurowski\cite{Nurowski}, Agrachev and Zelenko\cite{AZ}\cite{Zelenko2}. 
For example, Cartan distributions are related to the problem of \lq\lq rolling balls" (Agrachev\cite{Agrachev2}, 
Bor and Montgomery\cite{BM}, Baez and Huerta\cite{BaH}). 

Let $D \subset TY$ be a Cartan distribution on $Y$. 
In \cite{IKY}, it is shown that abnormal geodesics (or abnormal extremals, singular curves) of $D$ 
form another five dimensional space $X$ with a cone structure $C \subset TX$, 
and if the cone structure is regarded as a control system, 
the space of \lq\lq abnormal geodesics" of the cone structure is naturally identified with the original space $Y$. 
In fact, the duality is studied via the Cartan prolongation $(Z, E)$ of the given Cartan distribution $(Y, D)$ 
and the double fibration 
$$
(Y, D) \xleftarrow{\ \pi_Y \ } (Z, E) \xrightarrow{\ \pi_X\ } (X, C). 
$$
Then we obtain the decomposition $E = L \oplus K$ by integrable subbundles $L = \Ker(\pi_{Y*}), K = \Ker(\pi_{X*})$ of $E$. 
Note that abnormal geodesics are called singular paths in \cite{IKY}. 

For a space with Cartan distribution, the construction of the space of abnormal geodesics and  
the existence of natural double fibration was mentioned for 
the first time in unpublished lecture notes by Bryant \cite{Bryant}. 
Nurowski \cite{Nurowski} constructed indefinite conformal structures of signature $(3, 2)$ on such spaces. 
Then, as an alternative construction, 
the same cone structure as in \cite{IKY} was explicitly used in \cite{AZ} before \cite{IKY}. 
In fact, in \cite{AZ}, the cone structure was given by a foliation on the space $P((D^2)^\perp) \subset P(T^*Y)$, 
while in \cite{IKY} we describe it from an essentially same foliation in the space $P(D) \subset P(TY)$, 
which is canonically identified with $P(D^2)^\perp$. 

The construction of double fibrations and cone structures reminds us 
the general notion of pseudo-product structures in the sense of Tanaka (\cite{Tanaka}\cite{Yatsui}\cite{YY}, see also \cite{IM}). 

A {\it pseudo-product structure} on a manifold $M$ is a distribution $E \subset TM$ on a manifold $M$ with  
a decomposition $E = L \oplus K$ into integrable subbundles. Then locally we have a double fibration 
$$
P \xleftarrow{\pi_P} M \xrightarrow{\pi_N} N
$$ 
to the leaf space $P$ of $L$ and $N$ of $K$ respectively. 
Further we have the cone structures $K \hookrightarrow TM \xrightarrow{\pi_{P*}} TP$ on $P$ and 
$L  \hookrightarrow TM \xrightarrow{\pi_{N*}} TN$ on $N$ respectively. 
Thus we are naturally led to study the relation between abnormal geodesics on $P, N$ and on $M$. 
Our duality theorem (Theorem \ref{Duality}) 
treats the natural general problem in the case of Cartan distributions. 
Moreover we observe an asymmetry on this duality in terms of abnormal geodesics (Theorem \ref{Asymmetry}). 

The abnormal geodesics for the Cartan distribution with the maximal symmetry, namely with $G_2$-symmetry, 
was calculated in \cite{Kitagawa}. 
The double fibration $Y \leftarrow Z \to X$ with $G_2$-symmetry
and the canonical geometric structures $E \subset TZ, D \subset TY$ and $C \subset TX$ have been explicitly constructed 
in \cite{IMT}. 
Then in \cite{IKY} we determined abnormal geodesics, by direct calculations in $G_2$-case
for control systems ${\mathbb E} : E \hookrightarrow TZ \to Z$, ${\mathbb D} : D \hookrightarrow TY \to Y$ and 
${\mathbb C} = {\mathcal E}/\pi_X : E \to TX \to X$. 
Theorems \ref{Duality} and \ref{Asymmetry} show, from the viewpoint of geometric control theory 
the duality and asymmetry hold, not only for $G_2$-case, but also for general Cartan distributions. 

Note that double fibrations and cone structures arise also in works by Doubrov and Zelenko \cite{DZ1}\cite{DZ2}\cite{DZ3}. 
In Sato and Yamaguchi's paper \cite{SY}, a kind of geometric structure is naturally constructed on a space of singular curves (abnormal geodesics). 

In \S \ref{Control systems}, we recall the basic notions in geometric control theory, equivalence, admissible controls, controllability, 
etc., and, in particular, we make clear the notion of \lq\lq quotients" of control systems. 
In \S \ref{Abnormal geodesics}, we introduce the notion of abnormal geodesics, or singular paths, of a given control system and 
show the constrained Hamiltonian system describing abnormal geodesics. 
The duality and asymmetry on abnormal geodesics (singular paths) for Cartan distributions (\cite{IKY}) is explained in \S \ref{Cartan distributions}. 
Then the exposition is performed along the general notion of \lq\lq pseudo-product structures" in \S \ref{Pseudo-product structures and control systems}. 
In fact we show that a pseudo-product structure induces locally but naturally a double fibration and cone-structures, and 
we discuss the controllability of cone-structures briefly. 
After recalling the optimal control problems and associated {\it normal} and {\it abnormal extremals} in \S \ref{Optimal control problems}, 
we give the constrained Hamiltonian systems describing both normal and abnormal extremals (geodesics) for pseudo-product structures with 
sub-Riemannian metrics and for associated cone-structures with metrics in \S \ref{Sub-Riemannian pseudo-product structures}. 
Distributions can be regarded as cone structures of special types. 
Then naturally we hope the basic constructions given in this paper 
can be helpful to treat cone structures arising from pseudo-product structures. 
In \S \ref{Sub-Riemannian Cartan distributions}, 
we apply the general theory presented in the previous sections to the case of Cartan distributions. 

We do not treat in this paper on the classification problem of double fibrations, cone-structures nor pseudo-product structures. 
Nevertheless we should mention that we naturally suppose that 
the exact classification result of Cartan distributions by Zhitomirskii \cite{Zhitomirskii} can be interpret 
in terms of the classification of double fibrations via Cartan prolongations which have canonical pseudo-product structures. 
Also we should mention that the notion of Cartan connections plays the central role in the classification problem of geometric structures. 
Morimoto \cite{Morimoto} constructed the canonical Cartan connections associated with sub-Riemannian structures. Then we naturally ask 
for the construction of Cartan connections on sub-Riemannian pseudo-product structures.

\

All manifolds and mappings are assumed to be of class $C^\infty$ unless otherwise stated. 
For an interval $I$, we say that an assertion holds for almost every $t \in I$ if 
it holds outside of some measure zero set of $I$. 

\section{Control systems}
\label{Control systems}

Let $M$ be a finite dimensional $C^\infty$ manifold. 
A {\it control system} on $M$ is given by 
a locally trivial fibration $\pi_{\mathcal U} : {\mathcal U} \to M$ over $M$ and 
a $C^\infty$ mapping $F : {\mathcal U} \to TM$ to the tangent bundle $TM$ such that the diagram
$$
\begin{array}{ccc}
{\mathcal{U}} \  & \ \xrightarrow{ \ F \ } &  \ TM 
\vspace{0.2truecm}
\\
\ \ \ \ {\mbox{\footnotesize {$\pi_{\mathcal U}$}}}\searrow\hspace{-0.2truecm} &  
& \hspace{-0.2truecm}\swarrow{\mbox{\footnotesize {$\pi_{TM}$}}}  \ \ \ \
\vspace{0.2truecm}
\\
    & M &  
\end{array}
$$
commutes, where $\pi_{TM} : TM \to M$ is the projection of tangent bundle (\cite{Agrachev}).  
Then the control system is simply written as
$$
{\mathcal U} \xrightarrow{F} TM \xrightarrow{\pi_{TM}} M. 
$$ 

If a local triviality ${\mathcal U}\vert_{V} \cong V \times U, V \subset M$ on $M$ is given, then 
the control system is given by a family of vector fields 
$\{ f_u\}_{u \in U}$, $f_u(z) := F(z, u)$, and the equation 
$\dot{z}(t) = f_{u(t)}(z(t))$. 

Let ${\mathbb C} : {\mathcal U} \xrightarrow{F} TM \xrightarrow{\pi_{TM}} M$ be a control system over $M$ 
and $V \subset M$ an open subset. Then the {\it restriction} ${\mathbb C}\vert_V$ of ${\mathbb C}$ to $V$ is the control system 
$$
{\mathcal U}\vert_V \xrightarrow{F} TV \xrightarrow{\pi_{TV}} V. 
$$

Let ${\mathbb C} : {\mathcal U} \xrightarrow{F} TM \xrightarrow{\pi_{TM}} M$ 
be a control system over a manifold $M$ and $\pi_N : M \to N$ a fibration. 
Suppose also the composition $\pi_N\circ \pi_{\mathcal U} : {\mathcal U} \to N$ is a fibration. 
Then we have a control system ${\mathbb C}/\pi_N$ over $N$ by 
$$
{\mathbb C}/\pi_N : {\mathcal U} \xrightarrow{\pi_{N*}\circ F} TN \xrightarrow{\pi_{TN}} N. 
$$
Note that $\pi_{TN}\circ (\pi_{N*}\circ F) = \pi_N\circ \pi_{TM}\circ F = \pi_N\circ \pi_{\mathcal U}$. 
We call the control system ${\mathbb C}/\pi_N$ the {\it quotient} of ${\mathbb C}$ by $\pi_N$. 

Any distribution $D \subset TM$, namely, a vector subbundle of the tangent bundle $TM$,  is regarded as a control system 
$$
{\mathbb D} \ : \ D \hookrightarrow TM \longrightarrow M, 
$$
naturally via the fibration $\pi = \pi_{TM}\vert_D$ and the inclusion $F : D \hookrightarrow TM$. 

Two control systems ${\mathbb C} : {\mathcal U} \xrightarrow{F} TM \xrightarrow{\pi_{TM}} M$ and 
${\mathbb C}' : {\mathcal U}' \xrightarrow{F'} TM' \xrightarrow{\pi_{TM'}} M'$ are called {\it equivalent} if 
the following diagram commutes
$$
\begin{array}{ccccc}
{\mathcal U} & \xrightarrow{\ F\ } & TM & \xrightarrow{\pi_{TM}} & M
\\
\psi\downarrow \ \ \ & & \varphi_*\!\downarrow\ \ \  & & \ \ \downarrow\varphi
\\
{\mathcal U}' & \xrightarrow{\ F'\ } & TM' & \xrightarrow{\pi_{TM'}} & M'
\end{array}
$$
for some diffeomorphisms $\psi$ and $\varphi$. 
The pair $(\psi, \varphi)$ of diffeomorphisms is called an {\it equivalence} of the control systems ${\mathbb C}$ and 
${\mathbb C}'$. 

In this paper we suppose the fibre of $\pi : {\mathcal U} \to M$ is an open subset of $\R^r$. 

Given a control system ${\mathbb C} : {\mathcal U} \xrightarrow{F} TM \xrightarrow{\pi_{TM}} M$, 
an $L^\infty$ (measurable, essentially bounded) 
mapping $c : [a, b] \to {\mathcal U}$, $c(t) = (x(t), u(t)), x(t) \in M, u(t) \in \pi_{\mathcal U}^{-1}(x(t))$, 
is called a ${\mathbb C}$-{\it admissible control} or simply an {\it admissible control} 
if the differential equation 
$$
\dot{x}(t) = F(x(t), u(t)) \quad ({\mbox{\rm a.e.}}\  t \in [a, b]). 
$$
is satisfied. Here we suppose that $x(t) = \pi_{\mathcal U}\circ c(t)$ is absolutely continuous so that 
the above equation has a meaning. Then the curve $\pi_{\mathcal U}\circ c : [a, b] \to M$ turns to be a Lipschitz curve. 
We call the curve $\pi_{\mathcal U}\circ c$ the ${\mathbb C}$-{\it trajectory} 
associated to the admissible control $c$. 
We use the term \lq\lq path"  for a smooth ($C^\infty$) 
immersive trajectory regarded up to parametrisation as in \cite{IKY}. 

If the fibration $\pi_{\mathcal U} : {\mathcal U} \to M$ is trivial with a fibre $U$, namely if 
$\pi_{\mathcal U} : {\mathcal U} = M \times U \to M$ is the projection to the first factor, 
then the trajectory $x(t)$ associated to an admissible control $c : [a, b] \to M\times U$, $c(t) = (x(t), u(t))$
is uniquely determined by $u(t)$ and the initial point $x(a) \in M$. 

A control system ${\mathbb C}$ on a manifold $M$ is called {\it controllable} (resp. {\it locally controllable} at $z_0 \in M$) 
if any points $z_0, z_1 \in M$ are joined by a ${\mathbb C}$-trajectory. 
(resp. if, for any open neighbourhood $V \subset M$ of $z_0$, 
there exist open neighbourhood $W \subset V$ of $z_0$ such that 
any pair of points $z_1, z_2 \in W$ are jointed by a ${\mathbb C}\vert_V$-trajectory). 

\bel
\label{quotient-controllability}
For the quotient ${\mathbb C}/\pi_N$ of a control system ${\mathbb C}$ over a 
manifold $M$ by a fibration $\pi_N : M \to N$, we have
\\
{\rm (1) }
If $z : I \to M$ is a ${\mathbb C}$-trajectory, then $\pi_N\circ z : I \to N$ is 
a ${\mathbb C}/\pi_N$-trajectory. 
\\
{\rm (2) } 
If ${\mathbb C}$ is controllable, 
then ${\mathbb C}/\pi_N$ is controllable. 
\\
{\rm (2) } 
If ${\mathbb C}$ is locally controllable at a point $z_0 \in M$, 
then ${\mathbb C}/\pi_N$ is locally controllable at $\pi_N(z_0) \in N$. 
\enl

\beP
(1) : 
Let $z : I \to M$ is a ${\mathbb C}$-trajectory and $c : I \to {\mathcal U}$ an associated admissible control 
with $z = \pi_{\mathcal U}\circ c$, $\dot{z}(t) = F(c(t))$. Then $\pi_N\circ z = \pi_N\circ \pi_{\mathcal U}\circ c = 
\pi_{TN}\circ \pi_{N*}\circ F\circ c$ is a ${\mathbb C}/\pi_N$-trajectory which satisfies $\frac{\pa}{\pa t}(\pi_N\circ z)(t) = (\pi_{N*}\circ F)(c(t))$, 
associated with the same $c$. 
\\
(2) : 
Let $x_0, x_1 \in N$. Take $z_0, z_1 \in M$ such that $\pi_N(z_0) = x_0, \pi_N(z_1) = x_1$. 
There exists a ${\mathbb C}$-trajectory $z : [a, b] \to M$ 
such that $z(a) = z_0, z(b) = z_1$. 
Then $x(t) = \pi_N(z(t))$ is a ${\mathbb C}/\pi_N$-trajectory with $x(a) = x_0, x(b) = x_1$. 
\\
(3) : 
Let $V \subset N$ be any open neighbourhood of $x_0 = \pi_N(z_0)$. 
Since ${\mathbb C}$ is locally controllable at a point $z_0$, 
there exists an open neighbourhood $W' \subset \pi_N^{-1}(V)$ of $z_0$ such that 
any points in $W'$ are connected by a ${\mathbb C}\vert_{\pi_N^{-1}(V)}$-trajectory. 
Set $W = \pi_N(W')$. Then $W$ is an open neighbourhood of $x_0$ with $W \subset V$. 
Take any points $x_1, x_2 \in W$. Take $z_1, z_2 \in W'$ with $\pi_N(z_1) = x_1, \pi_N(z_2) = x_2$. 
Then there exists a ${\mathbb C}\vert_{\pi_N^{-1}(V)}$-trajectory 
$z : [a, b] \to {\mathcal U}\vert_{V}$ with $z(a) = z_1, z(b) = z_2$. 
Then $x = \pi_N\circ z$ is a $({\mathbb C}/\pi_N)\vert_V$-trajectory connecting $x_1$ and $x_2$. 
\enP

\section{Abnormal geodesics}
\label{Abnormal geodesics}

Fix $a, b \in \R$ with $a < b$. Then 
the totality ${\mathcal C}$ 
of admissible controls $c : [a, b] \to {\mathcal U}$ 
with a given initial point $z_0 \in M$, namely, with the condition $\pi_{\mathcal U}\circ c(a) = z_0$ 
form a Banach manifold. 
The {\it endpoint mapping} 
${\mathrm{End}} : {\mathcal C} \to M$ is defined by 
$$
{\mathrm{End}}(c) := \pi_{\mathcal U}\circ c(b).  
$$
The control system ${\mathbb C}$ is controllable if the endpoint mapping is surjective, for any $z_0$. 
We are concerning with the infinitesimal behaviour of the endpoint mapping. 

An admissible control $c : [a, b] \to {\mathcal U}$ with the initial point $\pi_{\mathcal U}(c(a)) = z_0$  
is called a {\it singular control} or an {\it abnormal control}, 
if $c \in {\mathcal C}$ is a singular point of ${\mathrm{End}}$, namely if 
the differential ${\mathrm{End}}_* : T_c{\mathcal C} \to T_{{\mathrm{End}}(c)}M$ is not surjective (\cite{BC}).  
If $c$ is a singular control, then the trajectory $\pi_{\mathcal U}\circ c$ is called a 
{\it singular trajectory}, {\it singular curve} or an {\it abnormal extremal}. 

A control system ${\mathbb C} : {\mathcal U} \xrightarrow{F} TM \xrightarrow{\pi_{TM}} M$ is 
called {\it trivial} if $F : {\mathcal U} \to TM$ is the zero map. Then the trivial system is not controllable at all, provided $M = \pt$. 

Given a control system 
${\mathbb C} : {\mathcal U} \xrightarrow{F} TM \xrightarrow{\pi_{TM}} M$, 
we set 
$$
\begin{array}{rrl}
{\mathcal U} \times_MT^*M & := & \{ ((z, u), (z', p)) \in {\mathcal U}\times T^*M \mid 
z = z' \}, 
\vspace{0.2truecm}
\\
& \cong & 
\{ (z, p, u) \mid z \in M, \ u \in \pi_{\mathcal U}^{-1}(z), \ p \in T_z^*M \}, 
\end{array}
$$
that is the fibre product of $\pi_{\mathcal U} : {\mathcal U} \to M$ and $\pi_{T^*M} : T^*M \to M$ over $M$. 

An equivalence $(\psi, \varphi)$ of control systems 
${\mathbb C} : {\mathcal U} \xrightarrow{F} TM \xrightarrow{\pi_{TM}} M$ and 
${\mathbb C}' : {\mathcal U}' \xrightarrow{F'} TM' \xrightarrow{\pi_{TM'}} M'$ induces the 
diffeomorphism 
$$
\psi\times_M\varphi^{-1*}  : {\mathcal U} \times_MT^*M \longrightarrow {\mathcal U}' \times_{M'}T^*M'
$$
defined by 
$$
(\psi\times_M\varphi^{-1*})((z, u), (z, p)) := (\psi(z, u), (\varphi(z), \varphi^{-1*}(p))). 
$$
Note that $\pi_{\mathcal U'}(\psi(z, u)) = \varphi(x)$. Moreover note that there exists the commutative diagram: 
$$
\begin{array}{ccccccc}
{\mathcal U} & \xrightarrow{\ F\ } & TM & \xrightarrow{\pi_{TM}} & M & \xleftarrow{\pi_{T^*M}} & T^*M
\\
\psi\downarrow \ \ \ & & \varphi_*\!\downarrow\ \ \  & & \ \ \downarrow\varphi &  & \ \ \downarrow\varphi^{-1*}
\\
{\mathcal U}' & \xrightarrow{\ F'\ } & TM' & \xrightarrow{\pi_{TM'}} & M' & \xleftarrow{\pi_{T^*M'}} & T^*M'. 
\end{array}
$$

We define the {\it Hamiltonian function} 
$$
H = H_{\mathbb C} \ : \ {\mathcal U} \times_MT^*M \longrightarrow \R
$$ 
of the control system ${\mathbb C}$ by 
$$
H((z, u), (z, p)) := \langle p, F(z, u)\rangle, \quad ((z,u), (z,p)) \in {\mathcal U} \times_MT^*M. 
$$
We will write $H(z, p, u)$ instead of $H((z, u), (z, p))$ for brevity. 

If $(\psi, \varphi)$ is an equivalence of control systems 
${\mathbb C}$ and ${\mathbb C}'$, then we have
$$
H_{\mathbb C} = H_{\mathbb C'}  \circ (\psi\times_M\varphi^{-1*}). 
$$
In fact, 
\begin{equation*}
\begin{split}
H_{\mathbb C'}  \circ (\psi\times_M\varphi^{-1*})((z, u), (z, p)) 
& = H_{\mathbb C'}(\psi(z, u), (\varphi(z), \varphi^{-1*}(p))) 
\\
& = \langle \varphi^{-1*}(p), F'(\psi(z, u)) \rangle 
\\
& = \langle \varphi^{-1*}(p), (\varphi_*)(F(z, u)) \rangle 
\\
& = \langle p, F(z, u) \rangle = H_{\mathbb C}((z, u), (z, p))
\end{split}
\end{equation*}

An admissible control $c : I \to {\mathcal U}$, $c(t) = (z(t), u(t))$ of the control system ${\mathbb C}$ 
is an abnormal control or a singular control if and only if 
there exists an absolutely continuous curve $\beta : I \to {\mathcal U} \times_MT^*M$, 
$\beta(t) = ((z(t), u(t)), (z(t), p(t)))$
satisfying, for any local triviality, the following constrained Hamiltonian system
$$
\begin{cases}
\ \dot{z}_i(t)  =  \ \ \dfrac{\pa H}{\pa p_i}(z(t), p(t), u(t)), \quad (1 \leq i \leq m)
\vspace{0.2truecm}
\\
\ \dot{p}_i(t)  =  - \dfrac{\pa H}{\pa z_i}(z(t), p(t), u(t)), \quad (1 \leq i \leq m)
\vspace{0.2truecm}
\\
\ \dfrac{\pa H}{\pa u_j}(z(t), p(t), u(t)) = 0, \quad (1 \leq j \leq r), \qquad 
p(t) \not= 0, 
\end{cases}
$$
for almost every $t \in I$, where $m = \dim(M)$. 

Note that abnormal (singular) controls are defined as \lq\lq singular points" of 
the end-point mapping which is one of principal objects to study in control theory. 
The above characterisation for them is a consequence of the calculation on the differential
of the endpoint mapping (see for example \cite{AS} pp.306--307). 

A curve $\beta : I \to {\mathcal U} \times_MT^*M$ satisfying the above constrained Hamiltonian system 
is called an {\it abnormal bi-extremal} and its projection $z : I \to M$ to $M$ is called an {\it abnormal extremal}
or a {\it singular trajectory}. 

We call a $C^\infty$ immersed abnormal extremal 
considered  up to $C^\infty$ parametrisations 
an {\it abnormal geodesic} or a {\it singular path} for the control system ${\mathbb C}$. 

For an equivalence $(\psi, \varphi)$ of control systems ${\mathbb C}$ and ${\mathbb C'}$, 
the abnormal ${\mathbb C}$-geodesics are mapped to the abnormal ${\mathbb C'}$-geodesics 
by the diffeomorphism $\varphi$. The class of abnormal geodesics is one of important invariants on the 
equivalence classes of control systems.

\section{Cartan distributions}
\label{Cartan distributions}

Let $Y$ be a $5$-dimensional manifold and 
$D \subset TY$ a distribution of rank $2$. 
Then $D$ is called a {\it Cartan distribution} if it has growth $(2, 3, 5)$, namely, if 
$\rank({\mathcal D}^{(2)}) = 3$ and $\rank({\mathcal D}^{(3)}) = 5$, where ${\mathcal D}^{(2)} = 
{\mathcal D} + [{\mathcal D}, {\mathcal D}]$ and ${\mathcal D}^{(3)} = 
{\mathcal D}^2 + [{\mathcal D}, {\mathcal D}^2]$ (see \cite{BH}). 
Here ${\mathcal D}$ means the sheaf of sections to $D$ and $[\ ,\ ]$ means the Lie bracket. 

Let $D \subset TY$ be a Cartan distribution. 
Consider the control system 
$$
{\mathbb D} : D \hookrightarrow TY \xrightarrow{\pi_{TY}} Y. 
$$
Then it is known that, for any point $y$ of $Y$ and for any direction $\ell \subset D_y$,  
there exists an abnormal $D$-geodesic (an abnormal geodesic or a singular path for ${\mathbb D}$) 
through $y$ with the given direction $\ell$ (\cite{LS}\cite{Montgomery}). 

Thus the abnormal $D$-geodesics form another five dimensional manifold $X$. 

Let 
$$
Z = PD = (D - 0)/\R^{\times} = \{ (y, \ell) \mid y \in Y, \ell \subset D_y, \dim(\ell) = 1 \}
$$ 
be the space of tangential lines in $D$. Then we have $\dim(Z) = 6$. 
Then $Z$ is naturally foliated by the liftings of abnormal $D$-geodesics, and we have locally double fibrations: 
$$
Y \xleftarrow{\ \pi_Y\ } Z \xrightarrow{\ \pi_X\ } X. 
$$
Note that, for each $y \in Y$, we have $\pi_Y^{-1}(y) = P(D_y)$, and, for each $x \in X$, 
$\pi_X^{-1}(x)$ is the locus in $Z$ of the lifted abnormal $D$-geodesics corresponding to $x$. 

Let 
$E \subset TZ$ be the {\it Cartan prolongation} of $D \subset TY$: 
For each $(y, \ell) \in Z$, $\ell \subset T_yY$, we set $E_{(y, \ell)} := \pi_{Y*}^{-1}(\ell)$. 
Then it is known that $E$ is a distribution of rank $2$ with growth $(2, 3, 4, 5, 6)$. 
If we put $L = \Ker(\pi_{Y*}), K = \Ker(\pi_{X*})$, then we have a decomposition 
$E = L \oplus K$ by integrable subbundles. 
Note that the decomposition is intrinsically obtained from the distribution $D \subset TY$. 

A cone field $C \subset TX$ is defined by setting, for each $x \in X$, 
$$
C_x := \bigcup_{z \in \pi_X^{-1}(x)}\pi_{X*}(L_z) \subset T_xX. 
$$


Thus, so far, we have distributions $D, E, L, K$ and the cone field $C$ from the double fibrations: 

$$
\begin{array}{ccccc}
TY & \hspace{-1truecm}\xleftarrow{\quad \pi_{Y*}\quad } \hspace{-2truecm}  & TZ & \hspace{-2truecm}\xrightarrow{\quad \pi_{X*}\quad } \hspace{-1truecm} & TX 
\vspace{0.2truecm}
\\
\bigcup & & \bigcup & & \bigcup
\vspace{0.2truecm}
\\
D & \xleftarrow{\ \pi_{Y*}\vert_E\ } & E = L \oplus K & 
\xrightarrow{\ \pi_{X*}\vert_E\ } & C 
\vspace{0.2truecm}
\\
\downarrow & & \downarrow & & \downarrow 
\vspace{0.2truecm}
\\
Y & 
\xleftarrow{\ \quad \pi_Y\quad\  } 
\hspace{-0.5truecm} 
& Z & 
\hspace{-0.5truecm} 
\xrightarrow{\ \quad \pi_X\quad  \ } 
& X
\end{array}
\vspace{0.2truecm}
$$
We regard the cone field $C \subset TX$ as a control system over $X$: 
$$
{\mathbb C} : L  \xrightarrow{\pi_{X*}\vert_L} TX \xrightarrow{\pi_{TX}} X.
$$ 
Then we have the following result: 
\bet 
\label{Duality}
{\rm (Duality \cite{IKY})}
Abnormal geodesics, namely singular paths, of the control system ${\mathbb C}$ 
are given by $\pi_X$-images of $\pi_Y$-fibres. 

Therefore, for any $x \in X$ and for any direction $\ell \subset C_x$, there exists uniquely 
an abnormal geodesic for the control system ${\mathbb C}$
passing through $x$ with the direction $\ell$. 
The original space $Y$ is identified with the space of singular paths on the space $X$, 
the space of singular paths on $Y$. 
\ent

Note that abnormal geodesics (singular paths) of ${\mathbb D} : D \hookrightarrow TY \to Y$ 
are given by $\pi_Y$-images of $\pi_X$-fibres.  

\ber
{\rm 
A different kind of duality was studied in \cite{IMT0} 
between a contact structure on a $3$-dimensional space and an indefinite conformal structure of signature $(2, 1)$ 
on another $3$-dimensional space via an Engel structure. 
Nevertheless, the duality via singular geodesics does not hold in that case, since 
contact structures do not have any abnormal geodesic, while the cone structure in that case actually has abnormal geodesics induced by 
the Engel line field. 
}
\enr

We are naturally led to suppose that there exists a kind of symmetry on $\pi_Y$-fibres and $\pi_X$-fibres. 
However the double fibration $Y \xleftarrow{\pi_Y} Z \xrightarrow{\pi_X} X$ is asymmetric 
in the following sense. 

An abnormal extremal $x(t)$ for $E \subset TZ$ is called {\it regular abnormal} 
if it is associated with an abnormal bi-extremal 
$(x(t), p(t), u(t))$ such that $p(t) \in E^{(2)\perp} \setminus E^{(3)\perp} \subset T^*Z$. 
An abnormal extremal $x(t)$ for $E \subset TZ$ is called {\it totally irregular abnormal} 
if any associated abnormal bi-extremal satisfies that $p(t) \in E^{(3)\perp} \subset T^*Z$. 
Then we have 

\bet
\label{Asymmetry}
{\rm (Asymmetry \cite{IKY})} 
An abnormal extremal, i.e. a singular trajectory for $E \hookrightarrow TZ \to Z$ 
is either a {\rm (}parametrisation of{\rm)} $\pi_Y$-fibre or a {\rm (}parametrisation of{\rm)} $\pi_X$-fibre. 
Each $\pi_Y$-fibre is regular abnormal, while 
each $\pi_X$-fibre is totally irregular abnormal. 
\ent

The main idea of the proof of Theorem \ref{Duality} is show that any abnormal geodesic of the cone structure 
$$
{\mathbb C} : L  \xrightarrow{\pi_{X*}\vert_L} TX \xrightarrow{\pi_{TX}} X, 
$$ 
regarded as a control system on $X$, lifts to an abnormal geodesic of 
$$
{\mathbb E} : E  \hookrightarrow TZ \xrightarrow{\pi_{TZ}} Z, 
$$
via $\pi_X : Z \to X$. Then we apply the first half of Theorem \ref{Asymmetry} 
to show that the lifted abnormal geodesic is either a $\pi_X$-fibre or 
a $\pi_Y$-fibre and then we conclude that the lifted abnormal geodesic must be a $\pi_Y$-fibre. 

We are going to present some of the results in \cite{IKY} in the following sections under a wider framework. 

\section{Pseudo-product structures and control systems}
\label{Pseudo-product structures and control systems}

Recall that a pseudo-product structure on a manifold $M$ is a distribution $E \subset TM$ on a manifold $M$ with  
a decomposition $E = L \oplus K$ into integrable subbundles. Then locally we have a double fibration 
$$
P \xleftarrow{\pi_P} M \xrightarrow{\pi_N} N
$$ 
to the leaf space $P$ of $L$ and $N$ of $K$ respectively. 
Then $K = \Ker(\pi_{N*}), L = \Ker(\pi_{P*})$. 

\bee
Let $E \subset TZ$ be the Cartan prolongation of a Cartan distribution $D \subset TY$. Then we have the 
intrinsic pseudo-product structure $E = L \oplus K$ as was shown in \S \ref{Cartan distributions}. 
\ene

Suppose $\pi_N\circ\pi_E : E \to N, \pi_N\circ\pi_L : L \to N, 
\pi_P\circ\pi_E : E \to P, \pi_P\circ\pi_K : K \to P$ are all fibrations. 
Then we are led to consider the following {\it five} control systems: 
$$
{\mathbb E} : E \hookrightarrow TM \xrightarrow{\pi_{TM}} TM, 
$$
which is a control system over $M$, 
$$
{\mathbb E}/\pi_N : E \hookrightarrow TM \xrightarrow{\pi_{N*}} TN, 
\ \ {\mbox{\rm and,}} \ \ 
{\mathbb L}/\pi_N : L \hookrightarrow TM \xrightarrow{\pi_{N*}} TN, 
$$
which are control systems over $N$, 
$$
{\mathbb E}/\pi_P : E \hookrightarrow TM \xrightarrow{\pi_{P*}} TP,  
\ \ {\mbox{\rm and,}} \ \ 
{\mathbb K}/\pi_P : K \hookrightarrow TM \xrightarrow{\pi_{P*}} TP, 
$$
which are control systems over $P$. 

Note that control systems ${\mathbb L}/\pi_P$ and ${\mathbb K}/\pi_N$ are reduced to be trivial. 

\bel
Let $z_0 \in M$. 
Suppose $E = L \oplus K \subset TM$ is bracket generating near $z_0$. 
Then the control system ${\mathbb E}$ {\rm (}resp. ${\mathbb E}/\pi_N, {\mathbb L}/\pi_N, {\mathbb E}/\pi_P, {\mathbb K}/\pi_P${\rm)}
is locally controllable at $z_0$ {\rm (}resp. locally controllable at $\pi_N(z_0), \pi_N(z_0), \pi_P(z_0), \pi_P(z_0)${\rm)}. 
\enl

\beP
By Chow-Rashevskii theorem \cite{Montgomery}, ${\mathbb E}$ is locally controllable at $z_0$. 
Then, by Lemma \ref{quotient-controllability}, ${\mathbb E}/\pi_N$ (resp. ${\mathbb E}/\pi_P$) is locally controllable at 
$\pi_N(z_0)$ (resp. at $\pi_P(z_0)$). 
To show the local controllability of ${\mathbb L}/\pi_N$ at $\pi_N(z_0)$, 
take a sufficiently small neighbourhood $V \subset N$ of $\pi_N(z_0)$ such that 
$\pi_N\circ \pi_{E} : E \to N$ and $\pi_N : M \to V$ are trivial over $V$ respectively. 
Moreover take a sufficiently small neighbourhood $V' \subset M$ of $z_0$ such that 
$\pi_N : V' \to V$ is a fibration and $L, K$ are trivial over $V'$. 
Then there is an open neighbourhood $W \subset V$ of $\pi_N(z_0)$ such that 
any pair of points $x_1, x_2 \in W$ is joined by an ${\mathbb E}/\pi_N$-trajectory $x : [a, b] \to V$. 
We take an associated ${\mathbb E}/\pi_N$-admissible control $c : [a, b] \to E\vert_{V'} = L\vert_{V'}\oplus K\vert_{V'}$ to $x$. 
Write $c(t) = (z(t), \ell(t), k(t))$ with $z : [a, b] \to V', \ell(t) \in L_{z(t)}, k(t) \in K_{z(t)}$. Then solve 
the differential equation $\frac{d}{dt}{\widetilde{z}}(t) = \ell(t)$ on curves $\widetilde{z} : [a, b] \to V'$ with an initial point $z(a)$. 
Then $\frac{d}{dt}(\pi_N\circ \widetilde{z})(t) = \pi_{N*}(\ell(t)) = \dot{x}(t)$. Then, by the uniqueness of the solution, 
we have $\pi_N\circ \widetilde{z} = x$. Therefore 
$x$ is an $({\mathbb L}/\pi_N)\vert_V$-trajectory connecting $x_1$ and $x_2$. Hence ${\mathbb L}/\pi_N$ is locally 
controllable at $\pi_N(z_0)$. Similarly we have that ${\mathbb K}/\pi_P$ is locally 
controllable at $\pi_P(z_0)$. 
\enP

\section{Optimal control problems}
\label{Optimal control problems}

Let ${\mathbb C} : {\mathcal U} \xrightarrow{F} TN \xrightarrow{\pi_{TN}} N$ be a control system on a manifold $N$. 

Given a smooth function $e :  {\mathcal U} \to \R$, 
consider an {\it optimal control problem} to minimise a {\it cost} functional 
$$
C := \int_a^b e(c(t)) dt 
$$
on admissible controls $c : [a, b] \to  {\mathcal U}$ with a fixed initial point $c(a) = (x_0, u_0)$, 
and a fixed endpoint $c(b) = (x_1, u_1)$. 
We express simply by the pair $({\mathbb C}, e)$ the optimal control problem. 
We call $e$ the {\it cost function} of the optimal control problem. 

Consider two control systems 
$$
{\mathbb C} : {\mathcal U} \xrightarrow{F} TN \xrightarrow{\pi_{TN}} N
\quad {\mbox{\rm and}} \quad 
{\mathbb C}' : {\mathcal U}' \xrightarrow{F'} TN' \xrightarrow{\pi_{TN'}} N', 
$$
and two cost functions $e : {\mathcal U} \to \R, e' : {\mathcal U}' \to \R$ respectively. 
Then two optimal control problems $({\mathbb C}, e)$ and 
$({\mathbb C}', e')$ 
are called {\it equivalent} if 
the diagram 
$$
\begin{array}{ccccccc}
\R & \xleftarrow{\ e\ } & {\mathcal U} & \xrightarrow{\ F\ } & TM & \xrightarrow{\pi_{TM}} & M
\\
|| & & \psi\downarrow \ \ \ & & \varphi_*\!\downarrow\ \ \  & & \ \ \downarrow\varphi
\\
\R & \xleftarrow{\ e'\ } & {\mathcal U}' & \xrightarrow{\ F'\ } & TM' & \xrightarrow{\pi_{TM'}} & M'
\end{array}
$$
commutes for some diffeomorphisms $\psi$ and $\varphi$. 
The pair $(\psi, \varphi)$ of diffeomorphisms is called an {\it equivalence} of the optimal control 
problems $({\mathbb C}, e)$ and $({\mathbb C}', e')$. 

\

The {\it Hamiltonian function} 
$$
H = H_{({\mathbb C}, e)} : ( {\mathcal U} \times_N T^*N) \times \R \longrightarrow \R, 
$$
of the optimal control problem $({\mathbb C}, e)$ is defined by 
$$
H_{({\mathbb C}, e)}(x, p, u, p^0) := H_{\mathbb C}(x, p, u) + p^0 e(x, u) = \langle p, F(x, u) \rangle + p^0 e(x, u).
$$
Here ${\mathcal U}\times_N T^*N$ means the fibre product 
of $\pi_{\mathcal U} : {\mathcal U} \to N$ and $\pi_{T^*N} : T^*N \to N$, namely, 
$$
{\mathcal U}\times_N T^*N := \{ (x, p, u) \mid (x, u) \in {\mathcal U}, (x, p) \in T^*N \}, 
$$
and $H_{\mathbb C}$ is the Hamiltonian function of the control system ${\mathbb C}$ as in \S \ref{Abnormal geodesics}. 

By the Pontryagin maximum principle, any solution (which is called a {\it optimal control}, a {\it minimal control} or {\it minimiser}) $c(t) = (x(t), u(t))$ of the optimal control 
problem is associated with a {\it bi-extremal} $(x(t), p(t), u(t), p^0)$ which 
is a solution of the 
constrained Hamiltonian system, for $H = H_{({\mathbb C}, e)}$, 
$$
\left\{
\begin{array}{l}
\dot{x}_i(t) = \dfrac{\pa H}{\pa p_i}(x(t), p(t), u(t), p^0), \ (1 \leq i \leq m), 
\vspace{0.3truecm}
\\
\dot{p}_i(t) = - \dfrac{\pa H}{\pa x_i}(x(t), p(t), u(t), p^0), \ (1 \leq i \leq m), 
\vspace{0.3truecm}
\\
\dfrac{\pa H}{\pa u_j}(x(t), p(t), u(t), p^0) = 0, \ (1 \leq j \leq r), \quad (p(t), p^0) \not= 0, \ p^0 \leq 0. 
\end{array}
\right.
$$
See \cite{PBGM}\cite{AS}. 

A bi-extremal is called {\it normal} (resp. {\it abnormal}) if $p^0 < 0$ (resp. $p^0 = 0$). 

A curve $x : I \to N$ (resp. a control $c = (x, u) : I \to {\mathcal U}$) is 
called a {\it normal extremal} (resp. an {\it abnormal extremal}) 
if it possesses a normal bi-extremal (resp. an abnormal bi-extremal) lift $(x(t), p(t), u(t), p^0)$. 
If a normal extremal (resp. an abnormal extremal) is regarded up to $C^\infty$ parametrisations, then it is called 
a {\it normal geodesic} (resp. an {\it abnormal geodesic}). Therefore we do not mean by a \lq\lq geodesic" 
an optimum but just a \lq\lq stationary point" of the cost functional in general. 

For an equivalence $(\psi, \varphi)$ of optimal control problems $({\mathbb C}, e)$ and $({\mathbb C'}, e')$, 
the ${\mathbb C}$-geodesics are mapped to the ${\mathbb C'}$-geodesics 
by the diffeomorphism $\varphi$. Thus class of normal geodesics is a natural invariant on the 
equivalence of optimal control problems. However it is difficult to describe all normal geodesics in general. 

\bee
{\rm
(\lq\lq generalised" sub-Riemannian geodesics). 
Let 
$X_1, \dots, X_r$ be vector fields over a manifold $N$. 
Consider the control system 
$$
\dot{x} = \sum_{i=1}^r u_i X_i(x). 
$$
Here $\pi : {\mathcal U} = N \times \R^r \to N$ is given by $\pi(x, u) = x$, 
and $F : {\mathcal U} \to TN$ is given by $F(x, u) = \sum_{i=1}^r u_i X_i(x)$. 
Moreover consider the optimal control problem to minimise 
the {\it energy} functional
$$
e = \int_I 
{\textstyle 
\frac{1}{2}\sum_{i=1}^r u_i(t)^2
}
\  dt. 
$$
It is known that the problem is equivalent to minimising the {\it length}: 
$$
\ell = \int_I \sqrt{
{\textstyle 
\sum_{i=1}^r u_i(t)^2
}
} \ dt. 
$$
If $X_1, \dots, X_r$ are linearly independent everywhere, then 
the problem is exactly to minimise the Carnot-Carath\'{e}odory distances in sub-Riemannian geometry (\cite{Montgomery}\cite{Bellaiche}). 

The Hamiltonian function of the optimal control problem is given by 
$$
H(x, p, u, p^0) = {\textstyle \sum_{i=1}^r \ \langle p, X_i(x)\rangle + \frac{1}{2} \, p^0\, (\sum_{i=1}^r u_i^2)}.
$$
Then the constraint $\frac{\pa H}{\pa u} = 0$ is equivalent to that 
$$
p^0 u_ j = - \langle p, X_j(x) \rangle, \ (1 \leq j \leq r). 
$$
For a normal extremal, we have $p^0 < 0$. Then we have 
$$
u_i = - \frac{1}{p^0} \langle p, X_j(x) \rangle, \ (1 \leq j \leq r). 
$$
Then 
$$
H = - \frac{1}{2p^0} \sum_{i=1}^r \ \langle p, X_i(x)\rangle^2. 
$$
From the linearity of Hamiltonian function on $(p, p^0)$, we can normalise $p^0$ to $-1$, so that 
$$
H = \frac{1}{2} \sum_{i=1}^r \ \langle p, X_i(x)\rangle^2. 
$$
Thus our formulation coincides with the ordinary normal Hamiltonian formalism in sub-Riemannian geometry (\cite{VG}\cite{LS}\cite{Montgomery}). 

For an abnormal extremal, the equation reduces to 
$$
\left\{
\begin{array}{l}
\dot{x}_i(t) = \dfrac{\pa H}{\pa p_i}(x(t), p(t), u(t)), \ (1 \leq i \leq m), 
\vspace{0.3truecm}
\\
\dot{p}_i(t) = - \dfrac{\pa H}{\pa x_i}(x(t), p(t), u(t)), \ (1 \leq i \leq m), 
\vspace{0.3truecm}
\\
\langle p, X_i(x)\rangle = 0, \ (1 \leq j \leq r), \quad p(t) \not= 0, 
\end{array}
\right.
$$
with $H(x, p, u) = H_{\mathbb C}(x, p, u) = \sum_{i=1}^r \ u_i\, \langle p, X_i(x)\rangle$. 
}
\ene

\section{Sub-Riemannian pseudo-product structures}
\label{Sub-Riemannian pseudo-product structures}

Let $E = L \oplus K \subset TM$ be a pseudo-product structure on $M$. 
Then we have associated with it, locally, a double fibration
$$
P \xleftarrow{\pi_P} M \xrightarrow{\pi_N} N. 
$$
Moreover we have associated five control systems to the pseudo-product structures in 
\S \ref{Pseudo-product structures and control systems}. 
We are going to treat several sub-Riemannian optimal control problems related to pseudo-product structures. 

Given Riemann metrics $g_L$ and $g_K$ on the integrable subbundles $L$ and $K$ respectively, 
we give the product metric $g_{E}$ of $g_L$ and $g_K$ on $E$: 
$$
g_E(v + w, v + w) := g_L(v, v) + g_K(w, w), \quad (v \in L_z, w \in K_z, z \in M). 
$$

\bee
{\rm 
Let $(E, g)$ be a sub-Riemannian structure on a distribution $E \subset TM$ and 
$L \subset E$ be an integrable subbundle of $E$. Then 
we set $K = L^{\perp}$ in $E$. If $K$ is integrable, 
then $E = L\oplus K$ and we have a sub-Riemannian pseudo-product structure. 
In particular, 
let $E$ be a sub-Riemannian Engel distribution and $L \subset E$ be the Engel line bundle (\cite{Montgomery}). 
Set $K := L^\perp \subset E$. Then $E = K \oplus L$ is a 
sub-Riemannian pseudo-product structure, $K, L$ being of rank one (integrable). 
}
\ene

%




\

Let $\dim M = m = k + \ell + q, \rank K = k$ and $\rank L = \ell$ so that $\rank E = r = k+\ell$. 
Then we put $p = k + q$ which is the dimension of $P$ 
and $n = \ell + q$ which is the dimension of $N$. 
Let us take local coordinates on $M$ in two ways: 
$$
y_1, \dots, y_{k+q}, v_1, \dots, v_{\ell}
{\mbox{\rm \quad and \quad}}
x_1, \dots, x_{\ell+q}, w_1, \dots, w_k, 
$$
where the projections $\pi_P, \pi_N$ are represented by 
$$
\pi_P(y_1, \dots, y_{\ell+q}, v_1, \dots, v_{\ell}) = (y_1, \dots, y_{k+q}), 
$$
$$
\pi_N(x_1, \dots, x_{\ell+q}, w_1, \dots, w_k) = (x_1, \dots, x_{\ell+q}), 
$$
respectively. 

We study the sub-Riemannian optimal control problems on the three control systems 
$$
{\mathbb E} : E \hookrightarrow TM \xrightarrow{\pi_{TM}} M
$$
over $M$, 
$$
{\mathbb E}/\pi_N : E \hookrightarrow TM \xrightarrow{\pi_{X*}} TN, 
\ {\mbox{\rm and,}} \ \
{\mathbb L}/\pi_N : L \hookrightarrow TM \xrightarrow{\pi_{X*}} TN
$$
over $N$. 
Optimal control problems on ${\mathbb E}/\pi_P$ and ${\mathbb L}/\pi_P$ are studied similarly. 

\

Note that $K = \Ker(\pi_{N*}) = \langle \frac{\pa}{\pa w_1}, \dots, \frac{\pa}{\pa w_k} \rangle$. 
Take local orthonormal frames $\xi_1, \dots, \xi_\ell$ of $L$ for given metric $g_L$ and $\eta_1, \dots, \eta_k$ of $K$ for given 
metric $g_K$ respectively: 
\begin{equation*}
\begin{split}
\xi_i(x, w) & \ = \ \sum_{1\leq j \leq n} c_{ij}(x, w)\frac{\pa}{\pa x_j} + \sum_{1 \leq \nu \leq k} e_{i\nu}(x, w)\frac{\pa}{\pa w_\nu},\quad (1 \leq i \leq \ell), 
\\
\eta_j(x, w) & \ = \ \sum_{1 \leq \nu \leq k} h_{j\nu}(x, w)\frac{\pa}{\pa w_\nu}, \quad (1 \leq j \leq k). 
\end{split}
\end{equation*}
Then we have a local orthonormal frame
$\xi_1, \dots, \xi_\ell, \eta_1, \dots, \eta_k$ of the sub-Riemannian metric on $E$ for the product metric $g_E$. 
Note that the $k\times k$-matrix $(h_{j\nu}(x, w))_{1 \leq j \leq k, 1 \leq \nu \leq k}$ is invertible. 

\

The system of canonical local coordinates of $T^*M$ is given by 
$$
(x, w; p, \psi) = (x_1, \dots, x_n, w_1, \dots, w_k; p_1, \dots, p_n, \psi_1, \dots, \psi_k). 
$$
For the control system ${\mathbb E}$, the Hamiltonian $H$ is given by 
$$
H_{\mathbb E}(x, w; p, \psi; a, b) = \sum_{1 \leq i \leq \ell} a_i H_{\xi_i}(x, w; p, \psi) 
+ \sum_{1 \leq i \leq k} b_iH_{\eta_j}(x, w; p, \psi), 
$$
where 
\begin{equation*}
\begin{split}
H_{\xi_i}(x, w; p, \psi)  & = \sum_{1\leq j \leq n} c_{ij}(x, w)p_j + \sum_{1\leq \nu \leq k} e_{i\nu}\psi_\nu, 
\\
H_{\eta_j}(x, w; p, \psi) & = \sum_{1\leq \nu \leq k} h_{j\nu}\psi_\nu, 
\end{split}
\end{equation*}
and $a = (a_1, \dots, a_\ell)$ (resp. $b = (b_1, \dots, b_k)$) are the fibre coordinates of $L$ (resp. $K$). 

We set 
$$
e_L := \frac{1}{2}\sum_{1\leq i\leq \ell} a_i^2, \quad e_K := \frac{1}{2}\sum_{1\leq j\leq k} b_j^2, 
$$
as the sub-Riemannian energy functions on $(L, g_L), (K, g_K)$ respectively and
$$
e_E : = e_L + e_K = \dfrac{1}{2}(\sum_{1\leq i\leq \ell} a_i^2 + \sum_{1\leq j\leq k} b_j^2), 
$$
as the total sub-Riemannian energy function on $(E, g_E)$. 

Then the Hamiltonian $H$ of the sub-Riemannian optimal control problem $({\mathbb E}, e_E)$ is given by 
$$
H := \sum_{1\leq i\leq \ell, 1 \leq j\leq n} a_ic_{ij}(x, w)p_j + \sum_{1\leq i\leq \ell, 1\leq\nu\leq k} 
a_ie_{i\nu}\psi_\nu + \sum_{1\leq j\leq k, 1\leq\nu\leq k} b_jh_{j\nu}\psi_\nu
+ \dfrac{1}{2}p^0(\sum_{1\leq i\leq \ell} a_i^2 + \sum_{1\leq j\leq k} b_j^2), 
$$
with variables $x, w, p, \psi$ and a non-positive parameter $p^0$. 

The constrained Hamiltonian system is given by 
\begin{equation*}
\begin{split}
\dot{x}_j & = \sum_{1\leq i \leq \ell} a_i c_{ij}, \quad (1 \leq j \leq n), 
\\
\dot{w}_j & = \sum_{1\leq i\leq \ell} a_ie_{i\nu} + 
\sum_{1\leq j\leq k} b_jh_{j\nu}, \quad (1 \leq j \leq k), 
\\
\dot{p}_{\nu} & = \ - \sum_{1\leq i \leq \ell, 1 \leq j \leq n}a_i\dfrac{\pa c_{ij}}{\pa x_{\nu}} p_j 
- \sum_{1\leq i \leq \ell, 1 \leq j \leq k}a_i\dfrac{\pa c_{ij}}{\pa x_{\nu}} \psi_j 
- \sum_{1\leq i \leq \ell, 1 \leq \mu \leq k}b_j\dfrac{\pa h_{i\mu}}{\pa x_{\nu}} \psi_\mu , 
\quad (1 \leq \nu \leq n), 
\\
\dot{\psi}_\nu & \ =  - \sum_{1\leq i \leq \ell, 1 \leq j \leq n}a_i\dfrac{\pa c_{ij}}{\pa x_{\nu}} p_j 
- \sum_{1\leq i \leq \ell, 1 \leq j \leq k}a_i\dfrac{\pa c_{ij}}{\pa x_{\nu}} \psi_j 
- \sum_{1\leq i \leq \ell, 1 \leq \mu \leq k}b_j\dfrac{\pa h_{i\mu}}{\pa x_{\nu}} \psi_\mu, 
\quad
(1 \leq \nu \leq k), 
\end{split}
\end{equation*}
\begin{equation*}
\sum_{1\leq j \leq n} c_{ij}p_j + \sum_{1\leq \nu \leq k} e_{i\nu}\psi_\nu + p^0a_i = 0, \ (1 \leq i \leq \ell), \quad 
\sum_{1\leq \nu \leq k} h_{j\nu}\psi_\nu + p^0b_j = 0, \ (1 \leq j \leq k), 
\end{equation*}
with $(p(t), \psi(t), p^0) \not= 0$. 

\

Consider the case $p^0 \not= 0$, i.e. the normal case. 
From the linearity of the equation on $(p, \psi, p^0)$, we can normalise $p^0 = -1$: If the equation is solved by 
$(p(t), \psi(t), p^0)$, then it is solved by $(\frac{1}{\vert p^0\vert}p(t), \frac{1}{\vert p^0\vert}\psi(t), -1)$. 
Then we have 
$$
a_i = \sum_{1\leq j \leq n} c_{ij}p_j + \sum_{1\leq \nu \leq k} e_{i\nu}\psi_\nu, \quad 
b_j = \sum_{1\leq \nu \leq k} h_{j\nu}\psi_\nu. 
$$

\

The optimal control problem $({\mathbb E}/\pi, e_{E})$ over $N$ is described as follows: 
The system of local coordinates of $N$ (resp. $T^*N$, $E$) 
is given by $x = (x_1, \dots, x_n)$, (resp. $(x; p) = (x; p_1, \dots, p_n)$, $(x, w; a, b) 
= (x, w; a_1, \dots, a_{\ell}, b_1, \dots, b_k)$). Then the system of local coordinates of 
$E\times_N T^*N$ is given by $(x, w; p; a, b)$. Note that $\dim(E\times_N T^*N) = 2n + 2k + \ell$. 
In this case, the mapping $F : E \to TN$ for the control system ${\mathbb E}/\pi$ is locally given by 
\begin{equation*}
\begin{split}
F(x, w; a, b)  & = \sum_{1 \leq i \leq \ell} a_i\, \pi_{X*}\xi_i(x, w) + 
\sum_{1 \leq j \leq k} b_j\, \pi_{N*}\eta_j(x, w) 
= \sum_{1 \leq i \leq \ell} a_i\, \pi_{N*}\xi_i(x, w), 
\end{split}
\end{equation*}
and the Hamiltonian function of the control system ${\mathbb E}/\pi$ and 
that for the optimal control problem $({\mathbb E}/\pi, e_E)$ are given respectively by 
\begin{equation*}
\begin{split}
H_{{\mathbb E}/\pi}(x, w; p; a, b) 
& 
= \ \langle p, \ \sum_{1 \leq i \leq \ell} a_i\, \pi_{X*}\xi_i(x, w) \rangle 
= \ \sum_{1 \leq i \leq \ell} a_i \langle p, \pi_{X*}\xi_i(x, w)\rangle  
\\
& = \ \sum_{1\leq i \leq \ell, 1 \leq j \leq n} a_i\, c_{ij}(x, w) p_j, 
\\
&
\\
H_{({\mathbb E}/\pi, e_E)}(x, w; p; a, b) & 
= \  H_{{\mathbb E}/\pi}(x, w; p; a, b) + e_E(a, b)
\\
& = \ \sum_{1\leq i \leq \ell, 1 \leq j \leq n} a_i\, c_{ij}(x, w) p_j 
+ \dfrac{1}{2}p^0(\sum_{1 \leq i \leq \ell} a_i^2 + \sum_{1 \leq j \leq k} b_j^2). 
\end{split}
\end{equation*} 
Here the control parameters 
are given by $w, a, b$. 
Note that $H_{{\mathbb E}/\pi}$ is independent of $b$, while $H_{({\mathbb E}/\pi, e_E)}$ 
does depend on $b$. The constrained Hamiltonian system for extremals is given by 
\begin{equation*}
\begin{split}
\dot{x}_j & = \sum_{1\leq i \leq \ell} a_i c_{ij}(x, w), \quad (1 \leq j \leq n), 
\\
\dot{p}_{\nu} & = - \sum_{1\leq i \leq \ell, 1 \leq j \leq n} a_i\dfrac{\pa c_{ij}}{\pa x_{\nu}}(x, w) p_j, 
\quad (1 \leq \nu \leq n), 
\end{split}
\end{equation*}
\begin{gather*}
 \sum_{1\leq i \leq \ell, 1 \leq j \leq n} a_i\dfrac{\pa c_{ij}}{\pa w_{\nu}} p_j = 0, 
(1 \leq \nu \leq k). 
\\
\sum_{1\leq j \leq n} c_{ij}p_j + p^0 a_i = 0, (1 \leq i \leq \ell), \quad p^0 b_j = 0, (1 \leq j \leq k). 
\end{gather*}

After normalising $p^0 = -1$, the constraints implies that 
$$
a_i = \sum_{1\leq j \leq n} c_{ij}p_j, (1 \leq i \leq \ell), \quad b_j = 0, (1 \leq j \leq k). 
$$

Suppose $\pi_{X*} : L \setminus \{ 0\} \to TX$ is an embedding. 
Then a normal extremal of $({\mathbb E}/\pi, e_E)$ lifts to a normal extremal of $({\mathbb E}, e_E)$ which is 
contained in a $\pi_P$-fibre, by setting $\psi = 0$. 
Then the lifted geodesic for $({\mathbb E}, e_E)$ contained in a $\pi_P$-fibre is a geodesic for the metric $g_L$. 
Thus we have 
\bep
\label{cone-geodesics}
Suppose $\pi_{X*} : L \setminus \{ 0\} \to TX$ is an embedding. 
Then the normal geodesics for the optimal control problem $({\mathbb E}/\pi_N, e_E)$ 
are obtained as $\pi_N$-images of Riemannian geodesics for the Riemannian metric $g_L$ in $L$-leaves. 
\enp
Similarly we have 
\bep
\label{cone-geodesics2}
Suppose $\pi_{Y*} : K \setminus \{ 0\} \to TY$ is an embedding. 
Then the normal geodesics for the optimal control problem $({\mathbb E}/\pi_P, e_E)$ 
are obtained as $\pi_P$-images of Riemannian geodesics for the Riemannian metric $g_K$ in $K$-leaves. 
\enp

\

Next we consider the optimal control problem $({\mathbb E}/\pi_N, e_L)$. 
The Hamiltonian function is given by 
$$
H_{({\mathbb E}/\pi, e_L)}(x, w; p; a, b)  \ = \ \sum_{1\leq i \leq \ell, 1 \leq j \leq n} a_i\, c_{ij}(x, w) p_j 
+ \dfrac{1}{2}p^0(\sum_{1 \leq i \leq \ell} a_i^2). 
$$
Then we have the following constrained Hamiltonian system: 
\begin{equation*}
\begin{split}
\dot{x}_j & = \sum_{1\leq i \leq \ell} a_i c_{ij}(x, w), \quad (1 \leq j \leq n), 
\\
\dot{p}_{\nu} & = - \sum_{1\leq i \leq \ell, 1 \leq j \leq n} a_i\dfrac{\pa c_{ij}}{\pa x_{\nu}}(x, w) p_j, 
\quad (1 \leq \nu \leq n), 
\end{split}
\end{equation*}
\begin{gather*}
\sum_{1\leq i \leq \ell, 1 \leq j \leq n} a_i\dfrac{\pa c_{ij}}{\pa w_{\nu}} p_j = 0, 
(1 \leq \nu \leq k), \quad 
\sum_{1\leq j \leq n} c_{ij}p_j + p^0 a_i = 0, (1 \leq i \leq \ell). 
\end{gather*}

After normalisation $p^0 = -1$, we have $a_i = \sum_{1\leq j \leq n} c_{ij}p_j$. 

\

We consider the optimal control problem $({\mathbb L}/\pi_N, e_L)$. 
The Hamiltonian function is given by 
$$
H_{({\mathbb L}/\pi, e_L)}(x, w; p; a)  \ = \ \sum_{1\leq i \leq \ell, 1 \leq j \leq n} a_i\, c_{ij}(x, w) p_j 
+ \dfrac{1}{2}p^0(\sum_{1 \leq i \leq \ell} a_i^2). 
$$
Then we have the following constrained Hamiltonian system: 
\begin{equation*}
\begin{split}
\dot{x}_j & = \sum_{1\leq i \leq \ell} a_i c_{ij}(x, w), \quad (1 \leq j \leq n), 
\\
\dot{p}_{\nu} & = - \sum_{1\leq i \leq \ell, 1 \leq j \leq n} a_i\dfrac{\pa c_{ij}}{\pa x_{\nu}}(x, w) p_j, 
\quad (1 \leq \nu \leq n), 
\end{split}
\end{equation*}
\begin{gather*}
\sum_{1\leq i \leq \ell, 1 \leq j \leq n} a_i\dfrac{\pa c_{ij}}{\pa w_{\nu}} p_j = 0, 
(1 \leq \nu \leq k), \quad 
\sum_{1\leq j \leq n} c_{ij}p_j + p^0 a_i = 0, (1 \leq i \leq \ell)
\end{gather*}

After normalisation $p^0 = -1$, we have $a_i = \sum_{1\leq j \leq n} c_{ij}(x, w)p_j$. 

\

We observe there there is no relation between normal geodesics for $({\mathbb E}, e_E)$ and those for 
$({\mathbb E}/\pi_N, e_L)$ (resp. $({\mathbb L}/\pi_N, e_L)$, because the latter problem has no information on 
$g_K$.

\

Now consider abnormal geodesics $(p^0 = 0)$. 
Then we have: 

\bel
\label{lifting}
{\rm (\cite{IKY})} \ 
Let $\gamma : I \to N$ be an abnormal $({\mathbb E}/\pi_N)$-trajectory. Suppose that 
there exists a Lipschitz abnormal bi-extremal $\beta : I \to E\times_N T^*N$ corresponding to $\gamma$. 
Then there exists an abnormal ${\mathbb E}$-trajectory $\widetilde{\gamma} : 
I \to M$ such that $\pi\circ \widetilde{\gamma} = \gamma$. 
\enl

In Lemma \ref{lifting}, 
we pose the Lipschitz condition on abnormal bi-extremals, because we have to regard some of control parameters 
as state variables. 

\

\noindent
{\it Proof of Lemma \ref{lifting}}: 
Let $\beta(t) = (x(t), w(t); p(t); a(t), b(t))$ be an abnormal bi-extremal for ${\mathbb E}/\pi_N$. 
Note that $b(t)$ can be taken as arbitrary $L^\infty$ function. 
Suppose $\beta(t)$ is Lipschitz. Then in particular the control $w(t)$ is Lipschitz. 
Replace $\mu(t)$ by $\dot{w}(t)$, which is of class $L^\infty$, and take $\psi(t) = 0$.
We set $\widetilde{\beta}(t) = (x(t), w(t); p(t), 0; a(t), \dot{w}(t))$. Then $\widetilde{\beta}$ is 
an abnormal bi-extremal for ${\mathbb E}$ and $\widetilde{\gamma}(t) = (x(t), w(t))$ 
is a lift of $\gamma(t) = x(t), \pi\circ \widetilde{\gamma} = \gamma$. 
\QED

\

The projection $\pi_L : E = L \oplus K \to L$ induces the  
projection $\rho : E\times_N T^*N \to L \times_N T^*N$. Then we have 

\bel
\label{lifting2}
{\rm (\cite{IKY})} \ 
A curve $\beta : I \to E \times_N T^*N$ is an abnormal bi-extremal {\rm (}resp. 
a Lipschitz abnormal bi-extremal{\rm )} for ${\mathbb E}/\pi$ if and only if 
$\rho\circ\beta : I \to L \times_N T^*N$ is an abnormal bi-extremal {\em (}resp. 
a Lipschitz abnormal bi-extremal{\rm )} for ${\mathbb L}$. Moreover any abnormal bi-extremal
{\rm (}resp. a Lipschitz abnormal bi-extremal{\rm )} $\overline{\beta} : I \to L \times_N T^*N$ for ${\mathbb L}$ 
is written as $\rho\circ\beta$ by an abnormal bi-extremal {\rm (}resp. 
a Lipschitz abnormal bi-extremal{\rm )} $\beta : I \to E \times_N T^*N$ for ${\mathbb E}/\pi$. 
\enl

\Proof
The system of local coordinates of $L \times_N T^*N$ is given 
$(x, w; p ; a)$.  The constrained Hamiltonian system for ${\mathbb L}$ 
is of the same form for ${\mathbb E}/\pi$ if the in-efficient component $b$ is deleted. 

We have that $\beta(t) = (x(t), w(t); p(t); a(t), bt))$ is an abnormal bi-extremal (resp. 
a Lipschitz abnormal bi-extremal) for ${\mathbb E}/\pi$ if and only if 
$\rho\circ\beta(t) = (x(t), w(t); p(t); a(t))$ is an abnormal bi-extremal (resp. 
a Lipschitz abnormal bi-extremal) for $L$. Therefore we have the first assertion. 
The second assertion is clear. 
\QED

\

These Lemmata are applied also to the systems ${\mathbb E}/\pi_P$ and ${\mathbb K}/\pi_P$. 
Using these general results, we have shown our main Theorem \ref{Duality} in \cite{IKY}.

\section{Sub-Riemannian Cartan distributions}
\label{Sub-Riemannian Cartan distributions}

We apply the general theory presented in the previous section to the case of Cartan distributions. 

Let $D \subset TY$ be a Cartan distribution and $g_D$ be a Riemannian metric on $D$. Then we have 
a sub-Riemannian Cartan distribution $(D, g_D)$. 
Let $E \subset TZ$ be the Cartan prolongation of $D$. Then we have the intrinsic pseudo-product structure 
$E = L \oplus K$, where $L$ and $K$ are of rank one (\S \ref{Cartan distributions}). 
Moreover there are induced natural metrics $g_L$ on $L$ and $g_K$ on $K$. 
In fact, for each $z \in Z$, $\pi_{Y*} : T_zZ \to T_{\pi_Y(z)}Y$ induces the injective linear map 
$K_z \to D_z$. Therefore the metric on $D_z$ induces a metric on $K_z$. Thus we obtain the metric $g_K$ 
from $g_D$. The metric $g_L$ is obtained from the following general construction. 

Let $V$ be a finite dimensional metric vector space of positive definite, namely, a Hilbert space with $\dim(V) < \infty$. 
Then we naturally define a Riemannian metric on the projective space $PV$ as follows. 
Let 
$$
SV = \{ x \in V \mid \ \Vert x\Vert = 1\} \subset V, 
$$
the unit sphere, which is a submanifold of codimension one. 
Then we have $PV = SV/\!\sim$ by identifying the pair of antipodal points. 
Let $\pi : SV \to PV$ denote the double covering. 
Take $[x] \in PV$ and $u, v \in T_{[x]}PV$. Define
$$
g_{PV [x]}(u, v) := g_{SV x}(\widetilde{u}, \widetilde{v}), 
$$
for $u = \pi_*(\widetilde{u}), v = \pi_*(\widetilde{v})$, 
$\widetilde{u}, \widetilde{v} \in T_xSV$. 
Then the metric $g_{PV}$ is well-defined because 
$$
g_{SV (-x)}(- \widetilde{u}, - \widetilde{v}) = g_{SV x}(\widetilde{u}, \widetilde{v})
$$
holds. 

\

In our case, the fibre $(\pi_Y)^{-1}(y)$ of $\pi_Y : Z \to Y$ over $y \in Y$ is identified with $P(D_y)$. 
Therefore $(\pi_Y)^{-1}(y)$ possesses the natural Riemannian metric, and, for each point $z \in Z$, 
$L_z = T_z((\pi_Y)^{-1}(\pi_Y(z)))$ has the natural metric. 
We endow $E$ with the product metric $g_E$ of $g_L$ and $g_K$. 

Let $\eta_1, \eta_2$ be an orthonormal frame of $(D, g_D)$. 
Take a system of local coordinates $y_1, y_2, y_3, y_4, y_5$ on $Y$. Then 
we have a system of local coordinates $y_1, y_2, y_3, y_4, y_5, v$ on $Z$ such that 
$\pi_Y(y, v) = y$, $v$ is the arc-length parameter of the fibres of $\pi_Y$. 
Then $\xi = \dfrac{\pa}{\pa v}$ is a unit frame of $L$, 
$$
\eta = (\cos v)\eta_1 + (\sin v)\eta_2 + \rho\dfrac{\pa}{\pa v}
$$
is a unit frame of $K$, for some function $\rho$ on $Z$, and then $\xi, \eta$ form the orthonormal frame of $E$. 
Write 
$$
\eta_1 = \sum_{1\leq i \leq 5} c_{1i}(y)\dfrac{\pa}{\pa y_i}, \quad \eta_2 = \sum_{1\leq i \leq 5} c_{2i}(y)\dfrac{\pa}{\pa y_i}. 
$$
Then 
$$
\xi = \dfrac{\pa}{\pa v}, \quad \eta = \sum_{1\leq i \leq 5}\{ (\cos v)c_{1i}(y) + (\sin v)c_{2i}(y)\}\dfrac{\pa}{\pa y_i} + \rho(y, v)\dfrac{\pa}{\pa v}. 
$$
The Hamiltonian function of the optimal control problem $({\mathbb E}, e_E)$ is given by 
$$
H_{({\mathcal E}, e_E)}(y, v; q, \varphi; \lambda, \mu; q^0) = 
\lambda\varphi + \mu\sum_{1\leq i \leq 5} \{(\cos v)c_{1i}(y) + (\sin v)c_{2i}(y)\} q_i + \mu \rho
+ \dfrac{1}{2}q^0(\lambda^2 + \mu^2), 
$$
where $(y, v; q, \varphi)$ is the system of local coordinates of $T^*Z$, $\lambda, \mu$ are the control parameters and 
$q^0$ is the non-positive constant. 
The constrained Hamiltonian system for $({\mathbb E}, e_E)$-extremal is given by 
\begin{equation*}
\begin{split}
\dot{y}_i & = \mu\{ (\cos v)c_{1i}(y) + (\sin v)c_{2i}(y)\}, \quad (1 \leq i \leq 5), 
\\
\dot{v} & = \lambda
\\
\dot{q}_j & = - \mu \sum_{1\leq i \leq 5} \{ (\cos v)\dfrac{\pa c_{1i}}{\pa y_j}(y) + (\sin v)\dfrac{\pa c_{2i}}{\pa y_j}(y)\}q_i - \mu\dfrac{\pa \rho}{\pa y_j}
\\
\dot{\varphi} & = - \sum_{1\leq i \leq 5} \mu \{ (- \sin v)c_{1i}(y) + (\cos v)c_{2i}(y)\}q_i - \mu\dfrac{\pa \rho}{\pa v}
\end{split}
\end{equation*}
$$
\varphi + q^0\lambda = 0, \quad \sum_{1\leq i \leq 5} \{ (\cos v)c_{1i}(y) + (\sin v)c_{2i}(y)\}q_i + \rho + q^0\mu = 0. 
$$

\

The Hamiltonian function of the optimal control problem $({\mathbb E}/\pi_Y, e_E)$ is given by 
$$
H_{({\mathbb E}/\pi_Y, e_E)}(y, v; q; \lambda, \mu; q^0) = 
\mu\sum_{1\leq i \leq 5} \{(\cos v)c_{1i}(y) + (\sin v)c_{2i}(y)\} q_i + \mu \rho
+ \dfrac{1}{2}q^0(\lambda^2 + \mu^2), 
$$
where $v, \lambda, \mu$ are control parameters. 
The constrained Hamiltonian system for $({\mathbb E}/\pi_Y, e_E)$-extremals is given by 
\begin{equation*}
\begin{split}
\dot{y}_i & = \mu\{ (\cos v)c_{1i}(y) + (\sin v)c_{2i}(y)\}, \quad (1 \leq i \leq 5), 
\\
\dot{q}_j & = - \sum_{1\leq i \leq 5} \mu \{ (\cos v)\dfrac{\pa c_{1i}}{\pa y_j}(y) + (\sin v)\dfrac{\pa c_{2i}}{\pa y_j}(y)\}q_i - \mu\dfrac{\pa \rho}{\pa y_j}, 
\quad (1 \leq i \leq 5), 
\end{split}
\end{equation*}
\begin{equation*}
\begin{split}
& \mu\sum_{1\leq i \leq 5} \{(- \sin v)c_{1i}(y) + (\cos v)c_{2i}(y)\} q_i = 0
\\
& q^0\lambda = 0, \quad \sum_{1\leq i \leq 5} \{ (\cos v)c_{1i}(y) + (\sin v)c_{2i}(y)\}q_i + \rho + q^0\mu = 0. 
\end{split}
\end{equation*}

We have the Hamiltonian function of the optimal control problem $({\mathbb E}/\pi_Y, e_K)$ by 
$$
H_{({\mathbb E}/\pi_Y, e_K)}(y, v; q; \lambda, \mu; q^0) = 
\mu\sum_{1\leq i \leq 5} \{(\cos v)c_{1i}(y) + (\sin v)c_{2i}(y)\} q_i + \mu \rho
+ \dfrac{1}{2}q^0\mu^2, 
$$
and the constrained Hamiltonian system for $({\mathbb E}/\pi_Y, e_K)$-extremals by 
\begin{equation*}
\begin{split}
\dot{y}_i & = \mu\{ (\cos v)c_{1i}(y) + (\sin v)c_{2i}(y)\}, \quad (1 \leq i \leq 5), 
\\
\dot{q}_j & = - \mu \sum_{1\leq i \leq 5} \{ (\cos v)\dfrac{\pa c_{1i}}{\pa y_j}(y) + (\sin v)\dfrac{\pa c_{2i}}{\pa y_j}(y)\}q_i - \mu\dfrac{\pa \rho}{\pa y_j}, 
\quad (1 \leq i \leq 5), 
\end{split}
\end{equation*}
\begin{equation*}
\begin{split}
\mu\sum_{1\leq i \leq 5} \{(- \sin v)c_{1i}(y) + (\cos v)c_{2i}(y)\} q_i = 0, 
\  \sum_{1\leq i \leq 5} \{ (\cos v)c_{1i}(y) + (\sin v)c_{2i}(y)\}q_i + \rho + q^0\mu = 0. 
\end{split}
\end{equation*}

Moreover we have $H_{({\mathbb K}/\pi_Y, e_K)}(y, v; q; \mu; q^0) = H_{({\mathbb E}/\pi_Y, e_K)}(y, v; q; \lambda, \mu; q^0)$
and the constrained Hamiltonian equation for $({\mathbb K}/\pi_Y, e_K)$ which is of the same form that for $({\mathbb E}/\pi_Y, e_K)$. 

\

Take another system of local coordinates $x_1, x_2, x_3, x_4, x_5, w$ on $Z$ such that 
$\pi_X(x, w) = x$ and take a orthonormal frame of $E$: 
$$
\xi = \sum_{1\leq i \leq 5} c_i(x, w)\dfrac{\pa}{\pa x_j} + f(x, w)\dfrac{\pa}{\pa w}, \quad 
\eta = \dfrac{\pa}{\pa w}. 
$$
We have, in this coordinates, 
\begin{equation*}
\begin{split}
H_{({\mathbb E}, e_E)}(x, w; p, \psi; a, b; p^0) 
& 
= 
\sum_{1\leq j \leq 5} ac_j(x, w)p_j + a f(x, w)\psi + b \psi + \dfrac{1}{2}p^0(a^2 + b^2), 
\\
H_{({\mathbb E}/\pi_X, e_E)}(x, w; p; a, b) 
& 
= 
\sum_{1\leq j \leq 5} ac_j(x, w)p_j + \dfrac{1}{2}p^0(a^2 + b^2), 
\\
H_{({\mathbb E}/\pi_X, e_L)}(x, w; p; a, b) 
& 
= 
\sum_{1\leq j \leq 5} ac_j(x, w)p_j + \dfrac{1}{2}p^0a^2, 
\\
H_{({\mathbb L}/\pi_X, e_L)}(x, w; p; a) 
& 
= 
\sum_{1\leq j \leq 5} ac_j(x, w)p_j + \dfrac{1}{2}p^0a^2, 
\end{split}
\end{equation*}
the control parameters are $(a, b)$, $(w, a, b)$, $(w, a, b)$, $(w, a)$ respectively. 
Though we can write down the corresponding constrained Hamiltonian equations in all cases, we omit to write down them here.

\

\begin{flushleft}
Goo ISHIKAWA, \\
Department of Mathematics, Hokkaido University, 
Sapporo 060-0810, Japan. \\
e-mail : ishikawa@math.sci.hokudai.ac.jp \\

\

Yumiko KITAGAWA, \\
Oita National College of Technology, 
Oita 870-0152, Japan. \\
e-mail : kitagawa@oita-ct.ac.jp

\

Wataru YUKUNO, \\
Department of Mathematics, Hokkaido University, 
Sapporo 060-0810, Japan. \\
e-mail : yukuwata@math.sci.hokudai.ac.jp

\end{flushleft}

\end{document}